\documentclass[preprint,12pt]{elsarticle}

\usepackage{amsmath,amssymb}
\usepackage{latexsym, color}
\usepackage{epsfig}

\begin{document}

\newtheorem{thm}{Theorem}[section]
\newtheorem{cor}[thm]{Corollary}
\newtheorem{lem}[thm]{Lemma}
\newtheorem{prop}[thm]{Proposition}
\newtheorem{defn}[thm]{Definition}
\newtheorem{rem}[thm]{Remark}
\newtheorem{nota}[thm]{Notation}
\newtheorem{Ex}[thm]{Example}
\def\nm{\noalign{\medskip}}

\bibliographystyle{plain}


\newcommand{\ds}{\displaystyle}
\newcommand{\pf}{\medskip \noindent {\sl Proof}. ~ }
\newcommand{\p}{\partial}
\renewcommand{\a}{\alpha}
\newcommand{\z}{\zeta}
\newcommand{\pd}[2]{\frac {\p #1}{\p #2}}
\newcommand{\norm}[1]{\left\| #1 \right \|}
\newcommand{\dbar}{\overline \p}
\newcommand{\eqnref}[1]{(\ref {#1})}
\newcommand{\na}{\nabla}
\newcommand{\Om}{\Omega}
\newcommand{\ep}{\epsilon}
\newcommand{\tmu}{\widetilde \epsilon}
\newcommand{\vep}{\varepsilon}
\newcommand{\tlambda}{\widetilde \lambda}
\newcommand{\tnu}{\widetilde \nu}
\newcommand{\vp}{\varphi}
\newcommand{\RR}{\mathbb{R}}
\newcommand{\CC}{\mathbb{C}}
\newcommand{\NN}{\mathbb{N}}
\renewcommand{\div}{\mbox{div}~}
\newcommand{\bu}{{\bf u}}
\newcommand{\la}{\langle}
\newcommand{\ra}{\rangle}
\newcommand{\Scal}{\mathcal{S}}
\newcommand{\Lcal}{\mathcal{L}}
\newcommand{\Kcal}{\mathcal{K}}
\newcommand{\Dcal}{\mathcal{D}}
\newcommand{\tScal}{\widetilde{\mathcal{S}}}
\newcommand{\tKcal}{\widetilde{\mathcal{K}}}
\newcommand{\Pcal}{\mathcal{P}}
\newcommand{\Qcal}{\mathcal{Q}}
\newcommand{\id}{\mbox{Id}}
\newcommand{\be}{\begin{equation}}
\newcommand{\ee}{\end{equation}}
\begin{frontmatter}



\title{Interfacial behaviors of the backward forward diffusion convection equations}


\author{Lianzhang Bao}

\address{School of Mathematics, Jilin University, Changchun, Jilin 130012, China, and School of Mathematical Science, Zhejiang University, Hangzhou 310027, China (lzbao@jlu.edu.cn)}

\begin{abstract}
This paper is devoted to the interfacial behaviors of a class of backward forward diffusion convection equations. Under the assumption that the equations have classical solutions in one dimension, we prove that the backward region shrinks and the forward region expands with a positive rate.\end{abstract}


\end{frontmatter}


\section{Introduction}\label{}

In this paper, we study the interfacial behaviors for a class of one dimensional backward-forward diffusion convection equations of the form:
\begin{equation}\label{eq:6001}
 u_t=(\Phi(u_x))_x + \Psi(x,u_x), \quad (x,t)\in Q_T=\Omega\times(0,T),
\end{equation}
where $\Omega\subseteq \mathbb{R}$ is an open bounded domain and $T>0$. The diffusion flux function $\Phi = \Phi(s)$ is assumed to be of a backward-forward type, including the well-known Perona-Malik type \cite{PM1992}. More precisely, we assume that there exist numbers $\alpha < \beta$ such that:
\begin{eqnarray}
 &&\Phi(s)\in C^{2}(\mathbb{R}), \label{eq:6002}
 \\
 && \Phi'(s)> 0,\forall s\in (\alpha,\beta),\label{eq:6003}
 \\
 &&\Phi'(s)< 0 ,\quad \forall s\in (-\infty, \alpha)\cup (\beta, +\infty). \label{eq:6004}
\end{eqnarray}
The non-monotonicity of $\Phi$ leads to a backward-forward partial differential equation \eqref{eq:6001} of parabolic type. A typical example of the backward-forward parabolic equation was considered by H$\ddot{o}$llig \cite{H1983}.
Similar equations to \eqref{eq:6001} with different types of nonmonotone flux functions $\Phi$ can also be found in fluid mechanics \cite{BBDPU},  image processing \cite{GP2012}, \cite{KY1}, \cite{KY2}, \cite{PM1992}, and gradient flows associated with nonconvex energy functionals \cite{BFG2006},\cite{S1991}.

\par Mathematical biology also provide many backgrounds to the backward-forward parabolic equations. Anguige et al. \cite{AS}, Bao and Zhou \cite{BZ}, Horstmann et al. \cite{HPO}, and Turchin \cite{T1989} all constructed the following backward forward equations:
\begin{equation}\label{eq:6005}
 u_t = (D(u)u_x)_x,
\end{equation}
where $u$ is the population density, and the diffusion coefficient $D(u)$ have positive and negative values corresponding to different physical interpretations.

\par Backward-forward parabolic equations with low order terms are also wildly studied. Wang et al. \cite{WNY} investigated the properties of the Young measure solutions of the backward-forward convection-diffusion equation.
 Kim and Yan \cite{KY3} proved the existence of infinitely many lipschitz weak solutions for the one-dimensional backward-forward parabolic equation of the following form
 \begin{equation}\label{eq:6006}
 u_t = (\sigma(u_x))_x + b(x,t)u_x  + c(x,t)u + f(x,t), \quad (x,t)\in Q_T=\Omega\times (0,T),
\end{equation}
where the monotonicity condition
\begin{equation}
 ({\sigma}(\xi)-{\sigma}(\zeta))\cdot (\xi-\zeta)\geq 0,
\end{equation}
is violated for some $\xi,\zeta \in \mathbb{R}^n$.

\par Interfacial behavior is a very important property of the backward-forward parabolic equation, especially in image processing that this is related to the movement of the detected edge. In the following, we denote
\begin{equation}\label{eq:6007}
Q^+:=\{(x,t)\in Q_T \quad\mbox{with}\quad u_x\in (\alpha,\beta)\}
\end{equation}
the forward (subcritical) region of Equation \eqref{eq:6001},
\begin{equation}\label{eq:6008}
Q^-:=\{(x,t)\in Q_T \quad\mbox{with}\quad u_x\in \mathbb{R}\setminus [\alpha,\beta]\}
\end{equation}
the backward (supercritical) region and
\begin{equation}\label{eq:6009}
Q^0:=\{(x,t)\in Q_T \quad\mbox{with}\quad u_x=\alpha \quad\mbox{or}\quad \beta\}
\end{equation}
the degenerate region. Note in the degenerate region, Equation \eqref{eq:6001} is a first order equation.

\par Classical solutions of \eqref{eq:6001} without low order terms have extensively been investigated in the last decade. Kawohl and Kutev \cite{KK1998} proved that global-in-time classical solutions exist if the initial conditions are in the forward region and under the framework of the classical solution, they proved the shrinking property of the union of supercritical and degenerate regions as time increasing and the supercritical region exists for all $t>0$ under some structural assumptions that the solution exists in one-dimension. Furthermore, Ghisi and Gobbino \cite {GG2010} proved similar results for one dimensional case and expanding property of supercritical regions in high dimensional case. The interfacial behaviors of a special backward-forward diffusion convection equation are also considered in \cite{BH2016}.
\par This paper is organized as follows. In section 2 we first state a result that is related to the interfacial behaviors of the degenerate parabolic equations which will be used in our proof of the main theorem. Second, we state our results for the backward-forward diffusion convection equation and the related interfacial behaviors in one dimension.
\section{Main results}\label{}
For simplicity, we assume that $\Phi$ satisfies condition \eqref{eq:6002},\eqref{eq:6003},\eqref{eq:6004} with $\Phi''(\alpha)>0,\Phi''(\beta)<0$,
\begin{equation}\label{eq:6013}
  \Psi(x,y)\in C^{1}(\mathbb{R}^2)\quad\mbox{and}\quad \Psi_x(x,\alpha)=A<0, \Psi_x(x,\beta)= B< 0.
\end{equation}
In order to prove our main theorem, we first state a result in \cite{GG2010} for the interfacial behaviors involving degenerate parabolic equations.
\begin{lem}\label{thm1}
 Let $x_1\leq x_2\leq x_3\leq x_4,$ and let $c_0, c_1, T, K, C$ be positive real numbers. Let $g: (0,c_0)\rightarrow (0,+\infty)$ be a continuous function such that
\begin{equation}\label{eq:6014}
 \lim_{\sigma\rightarrow 0^+} \frac{g(\sigma)}{\sqrt{\sigma}}=K.
\end{equation}
Let $f: (x_1,x_4)\times (0,T)\times [-c_1,c_1]^2 \rightarrow \mathbb{R}$ be a function such that $f(x,t,0,0)= C$ uniformly in $(x,t)$, namely
\begin{equation}\label{eq:6015}
 \lim_{\sigma\rightarrow 0^+}\sup\{|f(x,t,p,q)|: (x,t,p,q)\in (x_1,x_4)\times (0,T)\times[-\sigma,\sigma]^2\} =C.
\end{equation}
Let $v:(x_1,x_4)\times [0,T)\rightarrow \mathbb{R}$ be a function such that
\\
(c1) $v$ is continuous in $(x_1,x_4)\times [0,T);$
\\
(c2) $v(x,t)\geq 0$ for every $(x,t)\in (x_1,x_4)\times [0,T)$;
\\
(c3) $v(x,0)>0$ for every $x\in (x_2,x_3)$;
\\
(c4) $v_x(x,t)$ exists for every $(x,t)\in (x_1,x_4)\times (0,T);$
\\
(c5) setting
\begin{eqnarray}
 \mathcal{D} &:=&\{(x,t)\in (x_1,x_4)\times (0,T): 0<v(x,t)<c_0\},\label{eq:6102}
\\
 \mathcal{H}&:=&\{(x,t)\in (a,b)\times[0,T): x_2-k_0t<x<x_3+k_0t\}\label{eq:6103}
\end{eqnarray}
with $k_0=K\sqrt{C}$, we have that $v\in C^{2,1}(\mathcal{D})$, and
\begin{equation}\label{eq:6016}
 v_t\geq g(v)\{v_{xx}+f(x,t,v,v_x)\}\quad \forall (x,t)\in \mathcal{D}.
\end{equation}
Then $v(x,t)>0$ for every $(x,t)\in\mathcal{H}$.
\end{lem}
\par  The main strategy to prove the theorem is the comparison principle \cite{GG2010}.
\par Because Equation \eqref{eq:6001} involves only first order derivatives, we state our result as follows.

\begin{thm}\label{thm2}
 Let $\Phi\in C^{2}(\mathbb{R})$ be a function satisfying \eqref{eq:6002},\eqref{eq:6003},\eqref{eq:6004} with $\Phi''(\alpha)>0,\Phi''(\beta)<0$  and let $\Psi(x,y)$ satisfy \eqref{eq:6013}. Let $a\leq a_1\leq b_1\leq b$, and $T>0$ be a real number. Let $u\in C^1((a,b)\times [0,T))$ be a function satisfying \eqref{eq:6001} in $(a,b)\times[0,T),$ and
\begin{eqnarray}
 & \alpha<u_x(x,0)<\beta\quad \forall \quad x\in (a_1,b_1),\label{eq:6112}
  \\
  &u_x(x,0) <\alpha,\quad\forall\quad x\in[a,a_1),\quad u_x(x,0) >\beta,\quad\forall\quad x\in(b_1,b].\label{eq:6113}
\end{eqnarray}
Let $k_0 :=\sqrt{2|A\cdot\Phi''(\alpha)|},k_1 :=\sqrt{2|B\cdot\Phi''(\beta)|}$, and denote
\begin{equation}\label{eq:6012}
 \mathcal{G}:=\{(x,t)\in (a,b)\times[0,T): a_1-k_0t<x<b_1+k_1t\}.
\end{equation}
Then $\alpha<u_x(x,t)<\beta$ for every $(x,t)\in \mathcal{G}$.
\end{thm}

\par \pf Because $\Phi$ is an increasing function on $(\alpha, \beta)$, we first consider interfacial behaviors at $u_x = \beta$ and then consider the behaviors of the other side. Let us consider any function $\eta\in C^1(\mathbb{R})$ which is nondecreasing and where $\eta(\sigma)=\Phi(\sigma)$ for every $\sigma\in [\frac{\alpha + \beta}{2},\beta], \eta(\sigma)=\Phi(\beta)$ for every $\sigma\geq \beta$, $\eta(\sigma)=\Phi(\frac{3\alpha + \beta}{4})$ for $\sigma\leq \frac{3\alpha + \beta}{4}$. Then the function $\eta$ and $\Phi$ are invertible as a function from $(\frac{\alpha + \beta}{2},\beta)$ to $(\Phi(\frac{\alpha + \beta}{2}),\Phi(\beta))$. We can therefore define $g:(0,\Phi(\beta)-\Phi(\frac{\alpha + \beta}{2}))\rightarrow \mathbb{R}$ by setting
\begin{equation*}
 g(\sigma):= \Phi'(\eta^{-1}(\Phi(\beta)-\sigma)),\quad \forall \sigma\in (0,\Phi(\beta)-\Phi(\frac{\alpha + \beta}{2})),
\end{equation*}
which is a well defined continuous function. Finally we define
\begin{equation}\label{eq:6017}
 v(x,t) := \Phi(\beta)-\eta(u_x(x,t)),\quad \forall (x,t)\in (a,b)\times [0,T).
\end{equation}
If the constructed functions $g$ and $v$ satisfy conditions in Lemma \ref{thm1} which means $v(x,t) > 0$ in the expending region, and this is equivalent to $\eta(u_x(x,t)) < \Phi(\beta)$. By the definition of $\eta(\sigma)$, we obtain $u_x(x,t) <\beta$ in the region.

\par $\bf{Properties}$ $\bf{of}$ $g$. From the definition of function $g$, we can see $g: (0,\Phi(\beta)-\Phi(\frac{\alpha + \beta}{2}))\rightarrow (0,+\infty)$ is a continuous function. So we can square both the denominator and numerator, also change variable $\tau:= \eta^{-1}(\Phi(\beta)-\sigma)$, we have
\begin{eqnarray*}
 \lim_{\sigma\rightarrow 0^+}\frac{[g(\sigma)]^2}{\sigma}&=&\lim_{\sigma\rightarrow 0^+}\frac{[\Phi'(\eta^{-1}(\Phi(\beta)-\sigma))]^2}{\sigma} =\lim_{\tau\rightarrow \beta^-}\frac{[\Phi'(\tau)]^2}{\Phi(\beta)-\Phi(\tau)}
 \\
 &=& \lim_{\tau\rightarrow \beta^-}\frac{2\Phi'(\tau)\Phi''(\tau)}{-\Phi'(\tau)}=-2\Phi''(\beta) = 2|\Phi''(\beta)|,
\end{eqnarray*}
which proves \eqref{eq:6014} with $K:=\sqrt{2|\Phi''(\beta)|}.$
\par Properties of the function $v$. Because $\eta(\sigma)\leq \Phi(\beta)$ for every $\sigma\in \mathbb{R}$ and condition \eqref{eq:6112}, we have that $v$ satisfies $ (c1), (c2), (c3).$ By using standard interior regularity theory of parabolic equation, we have that $u$ is of class $C^\infty$ except at $ u_x(x,t)=\alpha $ and $ u_x(x,t)=\beta$. So we only need to consider the existence of $v_x(x,t)$ at $u_x(x_0,t_0)=\beta$, and the existence of $v_x(x,t)$ is trivial when $u_x(x,t) = \alpha$ because $h(\sigma)$ is constant for $\sigma \leq\frac{3\alpha +\beta}{4}$. We claim $v_x(x_0,t_0)=0$. Otherwise, there exists a sequence $\delta_k\rightarrow 0$ such that
\begin{equation}\label{eq:6018}
 |\frac{v(x_0+\delta_k, t_0)-v(x_0,t_0)}{\delta_k}|\geq \nu>0,\quad \forall k\in \mathbb{N}.
\end{equation}
\par By considering the subsequences, we can always assume that either $u_x(x_0+\delta_k, t_0)>\beta$ or $0<u_x(x_0+\delta_k, t_0)<\beta$ for every $k\in \mathbb{N}$. If $u_x(x_0+\delta_k, t_0)>\beta$, then the left side of Equation \eqref{eq:6018} is 0, which is a contradiction. In the second case, the right side of Equation \eqref{eq:6018} can be rewritten as
\begin{equation}\label{eq:6019}
 -\frac{\eta(u_x(x_0+\delta_k, t_0))-\eta(u_x(x_0,t_0))}{\delta_k} = -\frac{\Phi(u_x(x_0+\delta_k, t_0))-\Phi(u_x(x_0,t_0))}{\delta_k}.
\end{equation}
When $\delta_k\rightarrow 0$, Equation \eqref{eq:6019} tends to $(\Phi'(u_x))_x(x_0,t_0)$ which is 0, because the function $\Phi(u_x(x,t_0))$ attains its maximum at $x=x_0.$ Which is a contradiction to \eqref{eq:6018}.

\par In order to prove condition $(c5)$, we need to consider similar set defined as $\mathcal{D}$ in \eqref{eq:6102}. By using standard interior regularity theory of parabolic equation, we can see $v$ is regular and satisfies
\begin{equation}\label{eq:6020}
 v_t=- \eta'(u_x)u_{xt} = -\Phi'(u_x)[(\Phi(u_x))_x+ \Psi(x,u_x)]_x.
\end{equation}
Let $\Psi(x,u_x)=\Psi(x,\eta^{-1}(\Phi(\beta)-v))):=G(x,v)$ with $\Phi'(u_x)= \Phi'(\eta^{-1}(\Phi(\beta)-v))=g(v)$ and
\begin{equation}\label{eq:6021}
 [(\Phi(u_x))_x+ G_(x,v)]_x = (\Phi(u_x))_{xx}+ G_{x}(x,v) +G_{v}(x,v)v_{x}.
\end{equation}

Plugging these terms into \eqref{eq:6020}, we obtain:
\begin{equation}\label{eq:6022}
 v_t=g(v)\{v_{xx}- G_x(x,v) -G_v(x,v)v_{x}\}.
\end{equation}
By using Lemma \ref{thm1}, we have $\lim_{(v,v_x)\rightarrow 0^+}[G_x(x,v) +G_v(x,v)v_{x}]=B$, and conclude that the subcritical region expands with the rate $k_1 =\sqrt{2|B\cdot\Phi''(\beta)|}$.

\par We can prove the interfacial behavior for $u_x =\alpha$ in a similar way. Let $\eta_1\in C^1(\mathbb{R})$ be any nondecreasing smooth function where $\eta_1(\sigma)=\Phi(\sigma)$ for every $\sigma\in [\alpha, \frac{\alpha + \beta}{2}], \eta_1(\sigma)=\Phi(\alpha)$ for every $\sigma\leq \alpha$, $\eta_1(\sigma)=\Phi(\frac{\alpha + 3\beta}{4})$ for $\sigma\geq \frac{\alpha +3 \beta}{4}$. Then the function $\eta_1$ and $\Phi$ are invertible as a function from $(\alpha,\frac{\alpha + \beta}{2})$ to $(\Phi(\alpha), \Phi(\frac{\alpha + \beta}{2}))$. We can therefore define $g_1:(0,\Phi(\frac{\alpha + \beta}{2}) - \Phi(\alpha))\rightarrow \mathbb{R}$ by setting
\begin{equation*}
 g_1(\sigma):= \Phi'(\eta_1^{-1}(\Phi(\alpha)+\sigma)),\quad \forall \sigma\in (0,\Phi(\frac{\alpha + \beta}{2}) - \Phi(\alpha)),
\end{equation*}
which is a well defined continuous function. Finally we define
\begin{equation}\label{eq:6023}
 v^1(x,t) := \eta_1(u_x(x,t))-\Phi(\alpha),\quad \forall (x,t)\in (a,b)\times [0,T).
\end{equation}
Similar to the properties of $g$ and $v$ as in the $u_x =\beta$ side, we can obtain the following:
\begin{eqnarray*}
 \lim_{\sigma\rightarrow 0^-}\frac{[g_1(\sigma)]^2}{\sigma}&=&\lim_{\sigma\rightarrow 0^-}\frac{[\Phi'(\eta_1^{-1}(\Phi(\alpha)+\sigma))]^2}{\sigma} =\lim_{\tau\rightarrow \alpha^+}\frac{[\Phi'(\tau)]^2}{\Phi(\tau)-\Phi(\alpha)}
 \\
 &=& \lim_{\tau\rightarrow \alpha^+}\frac{2\Phi'(\tau)\Phi''(\tau)}{\Phi'(\tau)}=2\Phi''(\alpha),
\end{eqnarray*}
which proves \eqref{eq:6014} with $K:=\sqrt{2\Phi''(\alpha)}.$ The function $v^1$ is regular and satisfies
\begin{equation}\label{eq:6024}
 v^1_{t}= \eta_1'(u_x)u_{xt} = \Phi'(u_x)[(\Phi(u_x))_x+ \Psi(x,u_x)]_x.
\end{equation}
Let $\Psi(x,u_x)=\Psi(x,\eta_1^{-1}(\Phi(\alpha)+v))):=G^1(x,v)$ with $\Phi'(u_x)= \Phi'(\eta_1^{-1}(\Phi'(\alpha)+v))=g_1(v)$ and
\begin{equation}\label{eq:6025}
 [(\Phi(u_x))_x+ G^1(x,v)]_x = (\Phi(u_x))_{xx}+ G^1_{x}(x,v) +G^1_{v}(x,v)v_{x}.
\end{equation}

Plugging these terms into \eqref{eq:6020}, we obtain:
\begin{equation}\label{eq:6026}
 v^1_{t}=g_1(v)\{v^1_{xx}+ G^1_{x}(x,v) +G^1_{v}(x,v)v_{x}\}.
\end{equation}
By using Lemma \ref{thm1} again, we have $\lim_{(v,v_x)\rightarrow 0^+}[G^1_x(x,v) +G^1_v(x,v)v_{x}]=A$, and conclude that the subcritical region expands with the rate $k_1 =\sqrt{2A\cdot\Phi''(\alpha)}$.
\par
From the above results, we obtain the result that the subcritical region expands with lower bound positive rates which are depending on the function values at the vanishing point $u_x =\alpha$ and $u_x = \beta$.
\par\rightline{$\Box$}
\begin{rem}
When we consider the time behaviors of the supercritical region, by the time reverse we can change the supercritical region to the subcritical region and the lower bound expending rate is obtained as above which mean the supercritical region shrinking with a positive rate. The exact expending rate of subcritical region and shrinking rate of supercritical region may be different, which create a possibility that the emergence of a positive measured connected degenerate region from one point $x_0$ such that $u_x(x_0,0) =\alpha$ or $u_x(x_1,0)=\beta$ as time increasing.
\end{rem}
\section*{Acknowledgments}
 The research of this work was supported by China Postdoctoral Science Foundation--183816.

\end{document}